# Priscilla Greenwood: Queen of Probability

I.V. Evstigneev and N.H. Bingham

**Abstract:** This article contains the introduction to the special volume of Stochastics dedicated to Priscilla Greenwood, her CV and her list of publications.

This Festschrift celebrates Priscilla Greenwood, a remarkable mathematician, a bright personality and a valued friend and co-author for most of the contributors to the volume. For all of us, she has always been Cindy – Cindy Greenwood, the Queen of Probability. We are happy to honor her outstanding achievements with this Festschrift. This volume contains papers on Probability, Statistics and their applications. But these subjects do not cover all Cindy's areas of interest and those fields where she left her trace. Her scientific biography, which is sketched below, is unique and remarkable, and it reflects Cindy's beautiful mind and her true love of science.

Cindy Greenwood (C.G.) got into probability from operations research. She was a research assistant in OR at MIT as a graduate student. Feeling the need for stochastics in order to make any real progress, she took a first (and only) stochastic processes course there in 1960 with Henry P. McKean who lectured on the book he was writing with Itô. Also in that class were young Ramesh Gangolli, Don Dawson, Bert Fristedt, and Ron Getoor. After that year, C.G. moved to U. Wisconsin, Madison and finished a PhD in '63 with Josh Chover. The thesis was on a new method of prediction to replace the Wiener-Hopf method when data is on a compact interval. The method involved Schwartz distributions.

After two years in Durham, NC, where she taught many graduate courses in the black college then called North Carolina College, and joined seminars at Duke and UNC, Cindy moved to the U. of British Columbia (UBC), Vancouver. She was going to be in the Math Department there for 34 years, including 7 years of leave. She turned to questions about sample path properties of stochastic processes, particularly convergence of variations of Brownian motion and the variations of stable processes. This led to a collaboration with Bert Fristedt in 1970–71, when they were both on leave at U. Wisconsin, on variations of Lévy processes. In the 70s Cindy worked on several problems involving Wiener-Hopf factorization, in the random walk setting and unrelated to her thesis. She thought that it should be possible to produce Lévy processes from up-going and down-going factors, just as Brownian motion is produced from positive and negative drifts via random evolution. But this turned out not to work because the corresponding operators do not commute. The papers "Wiener-Hopf methods..." and "On Prabhu's factorization..." came out, and also, later, in a thwarted effort to define a natural bivariate exponential distribution, the papers with Moshe Shaked, "Fluctuation of random walk in $R^d$ ..." and "Dual pairs of stopping times for random walk". Moshe was a postdoc at Stanford when Cindy visited there on leave in 1976.

Another topic in the 70s was stopping times, especially in connection with regular variation. A resulting invention was the Martintote, a process which does the same thing with the asymptotic property of a process distribution that a Martingale does with the expected value. A paper showing that random stopping (usually) preserves regular variation, with Itrel Monroe, underscored the point. The Stanford visit began a collaboration with Sid Resnick, involving the somewhat novel idea of proving limit and structure theorems for Lévy processes starting with weak convergence of the point processes of jumps and then using convergence of continuous functionals to move to a (summed) process result. Cindy and Sid shed light on the structure of bivariate stable processes using this.

In 1977, Jim Pitman and C.G. were both at Churchill College and the Statistical Laboratory in Cambridge. Jim was interested in David Williams' theory of path decompositions and Itô's excursion theory for Markov processes, while C.G. was looking for a deeper probabilistic explanation of the Wiener-Hopf factorization for Lévy processes. Conversations around these interests led to a natural interpretation of fluctuation identities for random walks and Lévy processes by splitting at the maximum. This collaboration continued to provide a construction of local time and Poisson point processes, including Itô's excursion process, from nested arrays of independent and identically distributed sequences of random variables. In the summer of 1977, Jim and Cindy were invited by Elja Arjas to lecture at the probability summer school in Finland. A subsequent visit of C.G. to Finland produced another point process paper, this with Elja, on competing risks. Elja and family then visited Vancouver for a research year in the late 70s, as did Olav Kallenberg and his son.

Jef Teugels and C.G., who had first met on a street corner in Toronto enroute to a S.P.&A. Meeting at York in 1974, arranged to be at Cambridge at the same time, and worked together on a topic they named Harmonic Renewal Measures. This also had its roots in what Nick Bingham called "fluctuation theory for random walk", and was again related to Wiener-Hopf factorizations. At the end, Edward Omey was also involved.

In the early 80s Ed Perkins joined the Math Department at UBC, bringing his love of nonstandard methods, which Cindy had worked on earlier with Reuben Hersh. Together they studied Brownian local time on square root boundaries and excursions from moving boundaries, initially using non-standard arguments. But they yielded to social pressure and published with standard proofs. Also in the first part of the 80s, C.G. spent three summers plus half a year at the Steklov Math Institute in Moscow. She wrote papers with Albert Shiryaev which led both in the direction of statistics for stochastic processes. In addition there appeared a monograph on contiguity of stochastic process measures in a semi-martingale setting. There was also work on boundary crossing by Lévy processes with Alexander Novikov. In '83, C.G. visited Brighton and started work with Charles Goldie which led to papers on convergence of set-indexed processes. This was the start of recurring work on random fields.

Gerard Hooghiemstra, of Delft, visited C.G. in Vancouver, bringing his family and a problem on the storage of solar energy. This led to the definition of a family of operators interpolating between summation (sum) and supremum (sup), which they called alpha-sun operators because n is between m and p, and because the sun was involved. The domain of attraction theory for these operators was developed, and yielded new infinitely divisible distributions. The project was finished on a return visit of C.G. to Delft in '88.

Ron Doney also spent a year 1980–81 in Vancouver with his family. The topic of study was ladder variables of random walks, and for C.G. represented another return to the Wiener-Hopf circle of ideas. These ideas led, several years later, to a joint paper exploring the situation where the ladder times and heights belong to a bivariate domain of attraction, and also played a role in Ron's later solution of the problem of Spitzer's condition.

In the summer and fall of '86 C.G., was in Moscow, and was a member of the Soviet organizing committee of the first world congress of the Bernoulli Society in Tashkent. She traveled around the Soviet Union on these several visits and lectured in Russian in Novosibirsk, Kiev, Tashkent, and Yerevan, as well as Leningrad and Moscow. At the Vilnius conference in '86, C.G. lectured in Russian and gave her own simultaneous translation to English. This was followed by a visit to Leningrad, where a monograph with Mikhail Nikulin on Chi-squared testing was begun, and C.G. met Ildar Ibragimov. This led to a visit of Ildar to Vancouver and joint work on the Bahadur concept of efficiency.

In '87, with the sponsorship of Robert Serfling, C.G. received an NSF visiting professorship at Johns Hopkins. Andrew Barbour visited there briefly, resulting in joint work on Stein's method for Poisson approximation of random fields. C.G. gave a course at Hopkins on how to understand and use the new book of Jacod and Shiryaev on semi-martingales.

A collaboration lasting decades began with Wolfgang Wefelmeyer, who was also at Hopkins in '87 because of Serfling, on asymptotically efficient estimation in Le Cam's sense in the context of semi-martingales. Finding a cool reception of statistics in this setting, they focused on efficient estimation for particular processes under particular conditions, including counting processes, Markov chains, multivariate point processes, semi-Markov processes, random fields, partially specified models in various contexts, estimation near critical points, misspecified models, Monte Carlo Markov chain estimators, this last also with Ian McKeague. There were many visits back and forth between Vancouver and Germany, joint journeys to many conferences to present this work, and many papers.

Igor Evstigneev visited Vancouver in 1990. He and C.G. wrote a monograph on splitting and extrema of Markov random fields, published in the AMS Memoirs series. This is apparently the only monograph in the literature where Markov and strong Markov properties with respect to partial orderings, as well as analogues of splitting times, are systematically examined. A very different study in the same period with Mina Ossiander at Oregon State U., on functional convergence of evolving random fields, involved metric entropy conditions. After a sad misunderstanding with a major journal, where Cindy and Mina gave up too easily, the paper was published in the proceedings of the Sheffield Symposium on Applied Probability in 1991.

Ildar Ibragimov visited Vancouver in the early 90s, and the topic of study was Bahadur's concept of asymptotic efficiency. Much asymptotic efficiency work continued with Wefelmeyer through the 90s. In 1996 the Wiley book with Nikulin on chi-squared testing appeared after a decade of writing and rewriting.

Cindy's last year of leave from UBC., 1996–7, was spent half at the Weierstrass Institute in Berlin, visiting Michael Nussbaum, working on the Le Cam notion of asymptotic equivalence of statistical experiments, and the spring at Kent University in Canterbury visiting Howell Tong, studying the theory of dynamical systems as well as threshold models. Cindy's UBC postdoctoral colleague, Jaiming Sun, went also for the Berlin visit. Two papers by Sun and C.G. on Gibbs measures and the Ising model subsequently appeared in the Journal of Statistical Physics. Cindy went with Howell from Canterbury to a meeting in Oslo organized by Niels Chr. Stenseth who was interested in how to estimate the degree of synchronization of geographically related populations. Back at U. Kent, Cindy invented a suitable measure of synchronization and a model. Later in Vancouver she showed these ideas to statistician Dan Haydon who was at UBC then. A series of papers resulted, combining theory with data analysis, one on the Canadian Lynx data, one on muskrat and mink, also with Stenseth, and a third on the cycling of grey-sided vole data from Hokkaido, also with Stenseth and T. Sato using somewhat different modeling.

In the mid-90s UBC's Peter Wall Institute for Advanced Studies gave UBC faculty the challenge of forming a new interdisciplinary research group to study an otherwise inaccessible problem. Math Department people complained that the requirements eliminated them from the competition. The carrot was half a million dollars to fund a three-year project. Because of the work with Wefelmeyer on estimation of a parameter near a critical point, C.G. was aware of the special challenges and the ubiquity of critical points. The idea of the Crisis Points proposal was to study stochastic models with critical points arising from problems in various disciplines. After C.G.'s application, named "Crisis Points and Models for Decision", won the first Peter Wall support, colleague Ed Perkins complained, "You can't study that topic. That topic is everything!" In fact the Crisis Points Group collaborated on a variety of problems, involved about eight UBC faculty from various departments, and supported several grad students, post-docs and visitors. The post-docs included Ursula Müller from Bremen and Marty Anderies, now at ASU. Dan Haydon and C.G. started the synchronization work. UBC epidemiologist Steve Marion and Anders Martin Löf visiting from Stockholm worked with C.G. on stochastic epidemic models. Lawrence Ward in Psychology at UBC and Wefelmeyer worked with C.G. on stochastic resonance. The novelty of the stochastic resonance work in the Crisis Points Group was use of the asymptotically stable quantity Fisher information as a measure of the information produced by adding noise of a subthreshold signal, an idea that came from the statistical work with Wefelmeyer. One of the Crisis Points visitors was Frank Moss, a physicist from St. Lewis considered the father of stochastic resonance. He, Ward, C.G. and others wrote a stochastic resonance paper using Fisher information on the paddle fish, which became one of the most cited papers of the group. Stochastic resonance papers with combinations of Müller, Wefelmeyer and C.G. continued to appear in the years following C.G.'s move to Arizona State University (ASU) in the summer of 2000. The choice of ASU as a new venue resulted partly from the Crisis Points visit to Vancouver of Peter Killeen, who is in the Psychology Department at ASU.

Sveinung Edland came to ASU in '01 from Norway to work with C.G. for a year of his PhD research. He liked the question about the ubiquity of the 1/f phenomenon, which had been occupying C.G. for some time. The term 1/f refers to the phenomenon that the power spectral density of a process is nearly of the form "f to a negative power, the power being near 1", on a long interval of frequencies bounded away from zero. Ward and C.G. had worked on this during the Crisis Point years, had made some initial progress, but had not gotten much beyond what was already known. Sveinung and Cindy decided to find out what collection of Markov chains has this property, expecting this class to be broad. They found sufficient conditions in terms of the eigenstructure of a Markov transition matrix, and several examples. Efforts to publish this work where probabilists might find and pursue it failed. It appears in the usual place for papers on 1/f.

In the summer of '02, C.G. talked on measuring stochastic resonance using Fisher information at a conference in Brno in the Czech Republic. Petr Lansky, from the Institute of Physiology in Prague, attended, and a collaboration with visits both ways has yielded four papers, at present count, on stochastic resonance and other problems about the information in neuron-firing data. In '07, C.G. worked again on 1/f, this time in Prague and with the involvement of Lansky and his then research student, Lubomir Kostal. This work is ongoing. A recent survey on 1/f by L. Ward and C.G. can be found on the Scholarpedia website.

After two years in Arizona, C.G. spent a semester at Stockholm University lecturing on statistical genomics and working with Ola Hössjer on estimation of loci of disease-related genes. This work continued at ASU. For a year she was in the biology "dry lab" of Sudhir Kumar working on the evaluation of particular mutations in disease-related genes. During the year with Kumar, C.G. acquired the title Research Professor at ASU.

The Department of Mathematics and Statistics at ASU is strong in dynamical systems theory and mathematical biology. Stochastic dynamics had earlier aroused the interest of C.G., particularly in connection with the work of Rachel Kuske, now at UBC. Carlos Castillo-Chavez arrived at ASU in '04, joining the mathematical biology group and bringing the Mathematical and Theoretical Biology Institute which had been involving minority graduate students in epidemic modeling and analysis at Cornell for the previous ten years. With the resulting influx of graduate students and post-docs, mentors were needed. Carlos easily convinced C.G. along with K.P. Hadeler from Tubingen to become involved. This re-enforced the turn of C.G.'s research toward stochastic modeling in biology. She offered a graduate course under this title, and worked with post-docs and several students.

Rachel Kuske, post-doc Luis Gordillo and C.G. used multi-scale analysis to show that the damped oscillations in a deterministic susceptible-infective-removed (SIR) model become, in a stochastic analogue, sustained oscillations with a particular structure: sinusoids modulated by Ornstein-Uhlenbeck processes. Two papers which had been started with Anders Martin-Löf and Steve Marion were completed with Gordillo, one about the bimodal character of the total size distribution of an SIR epidemic, the second showing that in a non-homogeneous epidemic model total epidemic size does not depend on the pattern of infectiousness of the disease which may vary in the population.

The sustained oscillations in SIR epidemic models are clearly an example of a phenomenon associated with any dynamical system which has, in a deterministic form, oscillations damping to a locally stable equilibrium. For example, C.G. is pursuing, with Peter Rowat in Neural Computing at UCSD and Luis Gordillo, now at Puerto Rico, a way to compute the inter-spike interval distribution in stochastic neuron-firing models. The method may also help to explain the oscillatory nature of many populations.

With PhD student D. Tello and experimental psychologist Eddie Castañeda, C.G. is developing a model which will combine equations describing chemical production of neurotransmitter, e.g. dopamine, with the electrical drive and release which is neuron "firing". Apparently these two aspects of neural communication have never been modeled as one system. PhD student A. Mubayi and C.G. with others have produced two papers which model the dynamics of socially driven alcohol use and connect the model with college drinking data. A third paper explores an associated stochastic model. Azra Panjwani's M.Sc. thesis on the timely topic of the interplay between security and privacy under advancing technology is being rewritten by C.G. and research associate Xiaohong Wang as a paper for an ethics and technology journal.

Genomics and proteomics have expanded the possibilities for modeling of evolution, in contexts varying from speciation to drug resistance. Developing and analyzing such stochastic dynamic models tops the current collection of C.G.'s enthusiasms.

Cindy's career has been an impressive demonstration of mathematical power, sustained intellectual energy, enduring love of certain main themes, and at the same time an unusual openness to new ideas and new areas. Dr Samuel Johnson famously said that every man should strive to be an ornament to his profession. Cindy Greenwood has long been an ornament to her profession. She has touched all those who work in her field, by her prolific output of books and papers, by her teaching and professional service, by her scholarly qualities and her personal example. She has touched those fortunate enough to have worked with her closely by her warmth, her engaging personality, her fine sense of humour and her friendship. It is a pleasure to salute Cindy Greenwood and her achievements to date, and to wish her happiness, good health, and many years of productive scientific work to come.

Nick Bingham
Igor Evstigneev

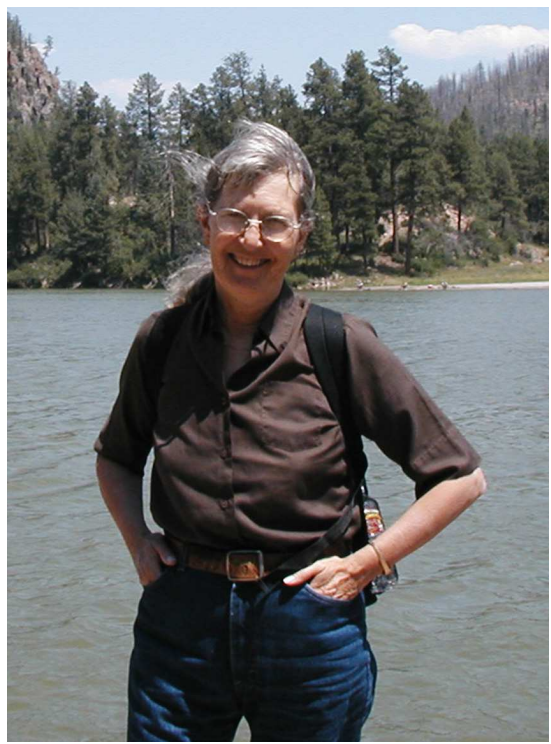

**Priscilla E. Greenwood**

Photo by Linda Gao

# Curriculum Vitae of
# Priscilla E. Greenwood

**Post-Secondary Education:**

| | | |
|---|---|---|
| Duke University, Durham, NC | BA | 1955–1959 |
| Massachusetts Institute of Technology | | 1959–1960 |
| University of Wisconsin-Madison | MA | 1961 |
| | PhD | 1963 |

**Employment:**

| | | |
|---|---|---|
| University of Wisconsin | Assistant Professor | 1963–1964 |
| North Carolina College, Durham | Associate Professor | 1964–1966 |
| University of British Columbia (UBC) | Professor (from 1985) | 1966–2000 |
| Arizona State University (ASU) Department of Mathematics and Statistics | Visiting Professor | Aug. 2000–2003 |
| Stockholm University | Visiting Professor | Sep.–Dec. 2002 |
| ASU Center for Evolutionary Functional Genomics | Senior Researcher | May 2003–May 2004 |
| ASU Department of Mathematics and Statistics | Research Professor | Jan. 2004–Present |

**Leaves of Absence from UBC and Visiting Positions in 1970–1997:**

| | | |
|---|---|---|
| University of Wisconsin | study leave | 1970–1971 |
| Mathematics Research Center, University of Wisconsin | Honorary Fellowship | Spring 1971 |
| Stanford University | study leave | July–Dec. 1976 |
| Churchill College, Cambridge, and Statistical Laboratory, Cambridge University, England | study leave | Jan.–July 1977 |
| University of California, Berkeley | study leave | Aug.–Dec. 1982 |
| University of Sussex, England | study leave | Jan.–May 1983 |
| Steklov Mathematical Institute, Moscow | visiting | July–Dec. 1986 |
| Johns Hopkins University | visiting | Jan.–Dec. 1987 |
| University of Delft, Netherlands | visiting | Jan.–June 1988 |
| University of Cologne and University of Heidelberg | study leave | July–Dec. 1990 |
| Colorado State University | study leave | Jan.–May 1991 |
| University of Heidelberg | study leave | May–Aug. 1991 |
| University of Stockholm | visiting | July.–Dec. 1993 |
| Weierstrass Institute, Berlin | study leave | Sep. 1996–Jan. 1997 |
| University of Kent, Canterbury, UK | study leave | Jan.–May 1997 |
| University of Siegen | study leave | June 1997 |

**Awards:**

Fellow, Institute of Mathematical Statistics, elected in August 1985.

Krieger-Nelson Prize for outstanding research record by woman mathematician from Canadian Mathematical Society, 2002.

Award of $500,000 from Peter Wall Institute of Advanced Studies for their first topic study: "Crisis Points and Models for Decision", 1997–2000.

**Teaching:**

Courses taught since 1966: many calculus, linear algebra, analysis, probability, and statistics courses.

Courses taught since 1999:

| | | |
|---|---|---|
| 1999/2000 | UBC | Probability |
| | UBC | Stochastic processes |
| 2000/2001 | ASU | Introduction to Statistics |
| | ASU | Probability |
| | ASU | Stochastic Processes |
| 2001/2002 | ASU | Calculus |
| | ASU | Probability |
| | ASU | Stochastic Processes |
| | ASU | Statistical Genomics |
| 2002/2003 | Stockholm University | Statistical Genomics |
| | ASU | Statistical Genomics Problems in Computational Bioinformatics |
| 2005/2006 | ASU | Stochastic Modeling in Biology Neurons and Neural Networks (co-taught) |
| 2006/2007 | ASU | Probability |

x

**Post-Docs Supervised (Since 1998):**

| | |
|---|---|
| Ursula Müller | Jan.–Aug. 1998 |
| John M. Andries | Sept. 1998–July 1999 |
| Greg Lewis | Sept. 1999–July 2000 |
| Luis Gordillo | 2004–2007 |

**PhD Students:**

| | |
|---|---|
| Suo Hong Chew | 1980 (co-supervised) |
| Harold Ship | 2004 (co-supervised) |

**MSc Students:**

| | |
|---|---|
| Glean Cooper | 1975 |
| David Szabo | 1987 |
| Brian Leroux | 1986 |
| Kongning Liu | 1989 |
| Harold Ship | 1993 |

# Priscilla E. Greenwood: Publications

1. A convolution equation on a compact interval. *Proc. Amer. Math. Soc.* **16** (1965) 8–13.
2. An asymptotic estimate of Brownian path variation. *Proc. Amer. Math. Soc.* **21** (1969) 134–138.
3. The variation of a stable path is stable. *Z. Wahrsch. Verw. Gebiete* **14** (1969) 140–148.
4. Variations of processes with stationary independent increments. *Z. Wahrsch. Verw. Gebiete* **23** (1972) 171–186 (with B. Fristedt).
5. Asymptotics of randomly stopped sequences with independent increments. *Ann. Probab.* **1** (1973) 317–321.
6. On Prabhu's factorization of Lévy generators. *Z. Wahrsch. Verw. Gebiete* **27** (1973) 75–77.
7. The Martintote. *Ann. Probab.* **2** (1974) 84–89.
8. Extreme time of processes with stationary independent increments. *Ann. Probab.* **3** (1975) 664–676.
9. Wiener-Hopf methods, decompositions, and factorisation identities for maxima and minima of homogeneous random processes. *Adv. Appl. Probab.* **7** (1975) 767–785.
10. Stochastic differentials and quasi-standard random variables. In: *Probabilistic Methods in Differential Equations (Proc. Conf., Univ. Victoria, Victoria, BC, 1974)*, Lecture Notes in Math., vol. 451, pp. 35–62. Berlin, Springer, 1975 (with R. Hersh).
11. Wiener-Hopf decomposition of random walks and Levy processes, *Z. Wahrsch. Verw. Gebiete* **34** (1976) 193–198.
12. Random stopping preserves regular variation of process distributions. *Annals of Probability* **5** (l977) 42–5l (with I. Monroe).
13. Fluctuations of random walk in $R^d$ and storage systems. *Adv. Appl. Probab.* **9** (1977) 566–587 (with M. Shaked).
14. Dual pairs of stopping times for random walk. *Ann. Probab.* **6** (1978) 644–650 (with M. Shaked).
15. A bivariate stable characterization and domains of attraction. *J. Multivariate Anal.* **9** (1979) 206–221 (with S. Resnick).
16. Fluctuation identities for Lévy processes and splitting at the maximum. *Adv. Appl. Probab.* **12** (1980) 893–902 (with J. Pitman).
17. Construction of local time and Poisson point processes from nested arrays. *J. London Math. Soc. (2)* **22** (1980) 182–192 (with J. Pitman).
18. Fluctuation identities for random walk by path decomposition at the maximum. *Adv. Appl. Probab.* **12** (1980) 291–293 (with J. Pitman).
19. Competing risks and independent minima: a marked point process approach. *Adv. Appl. Probab.* **13** (1981) 669–680 (with E. Arjas).
20. Point processes and system lifetimes. In: *Stochastic differential systems (Visegrád, 1980)*, Lecture Notes in Control and Information Sci., vol. 36, pp. 56–60. Berlin, Springer, 1981.
21. Harmonic renewal measures. *Z. Wahrsch. Verw. Gebiete* **59** (1982) 391–409 (with E. Omey and J.L. Teugels).